\documentclass[a4paper,12pt]{amsart}

\usepackage{amsmath,graphics}
\usepackage{amssymb}
\usepackage{amsfonts}
\usepackage{latexsym}
\usepackage{eucal}
\usepackage[dvips]{graphicx}
\usepackage{color}
\usepackage{enumerate}
\newtheorem{thm}{Theorem}[section]

\newtheorem{propo}[thm]{Proposition}
\newtheorem{lem}[thm]{Lemma}

\newtheorem{conj}[thm]{Conjecture}

\renewcommand{\Re}{{\rm Re}}
\renewcommand{\Im}{{\rm Im}}
\newcommand{\R}{\mathbb{R}}
\newcommand{\C}{\mathbb{C}}
\newcommand{\Z}{\mathbb{Z}}
\newcommand{\N}{\mathbb{N}}
\newcommand{\Q}{\mathbb{Q}}

\renewcommand{\H}{\mathbb{H}^2}

\newcommand{\U}{{\mathcal U}}
\newcommand{\half}{{\textstyle{\frac{1}{2}}}}

\newcommand{\D}{{\mathcal D}}
\newcommand{\A}{{\mathcal A}}

\begin{document}

\bibliographystyle{plain}
\title[Resonances of subgroups of arithmetic groups]{On the resonances of convex co-compact subgroups of arithmetic groups}
\author[Dmitry Jakobson]{Dmitry Jakobson}
\address{McGill university, Department of Mathematics and Statistics, Montreal, Quebec, Canada H3A2K6.}
\email{jakobson@math.mcgill.ca}
\author[Fr\'ed\'eric Naud]{Fr\'ed\'eric Naud}
\address{
Laboratoire d'Analyse non-lin\'eaire et G\'eom\'etrie (EA 2151)\\
Universit\'e d'Avignon et des pays de Vaucluse, F-84018 Avignon, France.
}
\email{frederic.naud@univ-avignon.fr}
\begin{abstract}
Let $\Lambda$ be a non-elementary convex co-compact fuchsian group which is a subgroup of an arithmetic fuchsian group.
We prove that the Laplace operator of the hyperbolic surface $X=\Lambda \backslash \H$ has infinitely
many resonances in an effective strip depending on the dimension of the limit set $\delta$. Applications to 
lower bounds for the hyperbolic lattice point counting problem are derived.
\end{abstract}

 \maketitle

\begin{section}{Introduction}
\noindent
Resonances are the natural replacement data for the missing eigenvalues of the Laplacian when dealing with non compact Riemannian
manifolds. In this paper we focus on hyperbolic surfaces.
Let $\H$ be the hyperbolic plane endowed with its standard metric of constant gaussian curvature $-1$.
Let $\Lambda$ be a geometrically finite discrete group of isometries acting on $\H$. This means
that $\Lambda$ admits a finite sided polygonal fundamental domain in $\H$. We will require that $\Lambda$ has no {\it elliptic} elements different from the identity and that the quotient $\Lambda \backslash \H$ is of {\it infinite hyperbolic area}. We assume in addition in this
paper that $\Lambda$ has no parabolic elements (no cusps).
Under these assumptions, the quotient space 
$X=\Lambda \backslash \H$ is a Riemann surface (called convex co-compact) whose geometry can be described as follows.
The surface $X$ can be decomposed into a compact surface $N$ with geodesic boundary, called the Nielsen region, on which $f$ ($f\geq 1$) infinite area ends $F_i$ are glued : the funnels. 
A funnel $F_i$ is a half cylinder 
$$F_i=(\R /l_i \Z)_\theta \times (\R^+)_t,$$ where $l_i>0$, with the warped metric 
$$ds^2=dt^2+\cosh^2(t)d\theta^2.$$
An important data related to the group $\Lambda$ is the Hausdorff dimension of its limit set $\delta(\Lambda)$, which in that case is also equal to the critical exponent of Poincar\'e  series.
Let $\Delta_X$ be the hyperbolic Laplacian on $X$. Its spectrum on $L^2(X)$ has been described completely by 
Lax and Phillips in \cite{LP2}. The half line $[1/4, +\infty)$ is the continuous spectrum, has no embedded eigenvalues.
The rest of the spectrum (point spectrum) is empty if $\delta\leq \half$, finite and starting at $\delta(1-\delta)$ if $\delta>\half$. The fact that the bottom of the spectrum is related to the dimension $\delta$ was first discovered by Patterson \cite{Patterson1} for convex co-compact groups.

\bigskip \noindent
By the preceding description of the spectrum, the resolvent 
$$R(s)=\left(\Delta_X-s(1-s) \right)^{-1}:L^2(X)\rightarrow L^2(X),$$
is therefore well defined and analytic on the half-plane $\{ \Re(s)> \half \}$ except at a possible finite set
of poles corresponding to the finite point spectrum. {\it Resonances} are then defined as poles of the meromorphic continuation of 
$$R(s):C_0^\infty(X)\rightarrow C^{\infty}(X)$$ 
to the whole complex plane. The set of poles is denoted by ${\mathcal R}_X$. This continuation is usually performed via the analytic Fredholm theorem
after the construction of an adequate parametrix. The first result of this kind in the more general setting
of asymptotically hyperbolic manifolds is due to Mazzeo and Melrose \cite{MazzMel}.
A more precise parametrix for surfaces was constructed by Guillop\'e and Zworski in \cite{GuiZwor,GuiZwor1}.
It should be mentioned at this point that in the infinite area case, resonances are spread all over the half plane $\{ \Re(s)< \half \}$,
in sharp contrast with the finite area case where resonances are known to be confined in a strip.
In this note, we adress the problem of finding resonances with the largest real part (which roughly speaking measures the decay
rate of resonant states). By Lax-Phillips theory, we know that when $\delta>\half$, all resonances (except the finite set of genuine $L^2$-eigenvalues) are in the half plane $\{ \Re(s)<\half\}$. On the other hand, if $\delta\leq \half$, a result of the second author \cite{Naud2}
shows that there exists $\epsilon>0$ such that ${\mathcal R}_X\cap \{\Re(s) \geq \delta-\epsilon\}=\{ \delta\}$. The constant $\epsilon$
obtained as a byproduct of a certain Dolgopyat estimate is barely effective. 

\bigskip
This dichotomy leads to the natural (and likely difficult) question: what is the actual size of the "essential spectral gap" i.e.
compute
$$ G(\Lambda):=\inf \left \{ \sigma <\delta\ :\  \{ \Re(s)\geq \sigma\}\cap {\mathcal R}_X\  \textrm{is finite} \right \}.$$
Clearly by the above discussion, we do have $G(\Lambda)\leq \half$ if $\delta>\half$ and $G(\Lambda)<\delta$ if $\delta\leq \half$.
We propose the following conjecture, see $\S 2$ for some heuristic justifications.
\begin{conj}
\label{deltaconj} 
Let $\Lambda$ be a convex co-compact group as above, then we have
$$G(\Lambda)=\frac{\delta}{2}.$$
\end{conj}
\noindent
This conjecture is consistent with the fact that in the finite volume 
case ($\delta=1$), a result of Selberg \cite{Selberg2} shows that 
there exist resonances in the strip $$\{\Re(s)<\half\}$$
 whose real parts are as close to $\half$ as we want. However, his 
strategy does not work in the infinite volume case, mostly because
 resonances dot not lie in a strip. In this paper, we prove the following.
 \begin{thm}
  \label{mainresult}
  Let $\Lambda$ be a convex co-compact Fuchsian group.
  \begin{itemize}
   \item If  $0<\delta\leq \half$, then we have $G(\Lambda)\geq 
\frac{\delta(1-2\delta)}{2}$.
   \item If  $\half< \delta$ and $\Lambda$ is a convex co-compact 
subgroup of an arithmetic group, then
    $G(\Lambda) \geq \frac{\delta}{2}-\frac{1}{4} $.
  \end{itemize}
 \end{thm}
 \noindent
In the elementary case, $\delta=0$, and this lower bound is known 
to be an equality since resonances can be explicitly computed, see
\cite{ChristZworski}. We also notice that as $\delta\rightarrow 0$, we 
have 
$\rho(\delta)=\frac{\delta}{2}+O(\delta^2)$ so our lower bound is, at 
least infinitesimally, in favor of the conjecture. The above estimate 
clearly beats previous results on strips with infinitely many resonances 
obtained for arbitrary convex co-compact hyperbolic manifolds in 
\cite{GuiNaud2}. If $\Lambda$ is non elementary and $\delta\neq \half$, 
then this proves that $G(\Lambda)>0$. 

It follows from a recent result of Lewis Bowen \cite{LBowen} (see 
also \cite{Gamburd, BK1} for the $\mathrm{SL}_2(\Z)$ case)
that in every co-finite or co-compact arithmetic Fuchsian group, 
one can find convex co-compact subgroup with dimension
arbitrarily close to $1$. So there are plenty of examples of 
convex co-compact surfaces with infinitely many resonances whose real
parts are arbitrarily close to $1/4$.

One of the motivations for this paper is the hyperbolic lattice point 
counting problem. Indeed,
a famous result of Lax-Phillips \cite{LP1} shows that for $\delta>\half$, 
for all $z,z' \in \H$, there exists a finite asymptotic expansion of the 
counting function (``${\rm d}$'' denotes hyperbolic distance)
$$N(T;z,z'):=\#\{\gamma \in \Lambda\ :\ {\rm d}(z,\gamma z')\leq T \},$$ 
which is as $T\rightarrow +\infty$
$$N(T;z,z')=\sum_j C_j(z,z')e^{\delta_jT}+O\left( T^{5/6} 
e^{(\delta+1)T/3}\right),$$
where $j$ runs over a finite set and each $\delta_j \in (\half,\delta]$ 
is related to the small eigenvalues of $\Delta_X$ through the formula
$\delta_j(1-\delta_j)=\lambda_j$, $\delta_0=\delta$. It is conjectured 
in the {\it finite volume } case that the optimal error term is
$O(e^{(1/2+\epsilon)T})$. In the infinite volume case, in view of our 
conjecture, the optimal error term should be $O(e^{(\delta/2+\epsilon)T})$ and 
one might expect a more sophisticated expansion involving resonances 
in addition to $L^2$ point spectrum.

As a byproduct of Theorem \ref{mainresult}, we obtain the following.
\begin{thm}
\label{lattice}
Let $\Lambda$ be a convex co-compact subgroup of an arithmetic group 
with $\delta>\half$. There exists a full measure subset 
${\mathcal G}\subset \H\times \H$ such that for all $(z,z')\in
 {\mathcal G}$ and all finite expansion of the form
$$\sum_j Q_j(T;z,z')e^{\delta_jT},$$
where $\delta_j \in \C$ and $Q_j(T;z,z') \in \C[T]$, we have for 
all $\epsilon>0$
$$\left \vert N(T;z,z')-\sum_j Q_j(T;z,z')e^{\delta_jT}\right 
\vert=\Omega\left( e^{(\delta/2-1/4-\epsilon)T}\right).$$
\end{thm}
\noindent
The notation $\Omega(\bullet)$ means, as usual in number theory, 
being not a $O(\bullet)$.
The plan of the paper is as follows. In section $\S2$, we motivate 
conjecture \ref{deltaconj} by some heuristics based
on thermodynamical formalism that might lead to a future proof.
In section $\S3$ we prove Theorem \ref{mainresult}. We first prove an 
approximate Trace formula
to take advantage of the known upper bounds on Selberg's zeta function. 
The result then follows from the prime geodesic theorem and
a mean square argument ``\`a la Selberg''. Arithmeticity is used when 
$\delta>\half$ to improve lower bounds due to exponentially growing 
multiplicities in the length spectrum. Section $\S 4$ is devoted to 
the application to the lattice point counting problem. Although it is 
a fairly direct
consequence of the previous result, some details have to be provided, 
in particular we explain in terms of the polar structure of the resolvent 
where the generic set of $(z,z')$ for which the lower bound holds comes from. 

\bigskip
Finally, we mention that using the trace formula from \cite{GuiNaud1} 
or the Patterson-Perry result on the divisor of the Selberg zeta 
function \cite{PatPer}, our result should extend straightforwardly 
to higher dimensional convex co-compact hyperbolic manifolds 
of Schottky type (the fractal Weyl upper bound from 
Guillop\'e-Lin-Zworski \cite{GLZ}, see estimate (\ref{fractal}) below, 
is critical in our proof). Similarly, if the fractal Weyl upper bound 
(\ref{fractal}) holds for surfaces with cusps, then our result should extend 
straightforwardly to that case.

\bigskip \noindent {\bf Acknowledgement.} We thank Peter Sarnak for 
pointing out reference \cite{LBowen}. Both Authors thank
for their hospitality the Banff center and McGill university where 
part of this work has been done. DJ is supported by NSERC, FQRNT
and Dawson fellowship. FN is supported by ANR METHCHAOS.
\end{section}
\begin{section}{Transfer operators and heuristics on the conjecture}
Recall that Selberg's zeta function $Z_\Lambda(s)$ is defined for
$\Re(s)$ large by the infinite product
$$Z_\Lambda(s):=\prod_{n\in \N}\prod_{\gamma \in {\mathcal P}}\left( 1-e^{-(s+n)l(\gamma)} \right),$$
where ${\mathcal P}$ is the set of prime closed geodesics on $X=\Lambda \backslash \H$ and if $\gamma \in {\mathcal P}$, 
$l(\gamma)$ is the length.  By a rough upper bound on the number of closed geodesics, convergence of the above product actually
occurs for all $\Re(s)>\delta$. Since $\Lambda$ is a convex co-compact group, it is known that the above product extends analytically
to an entire function on the whole complex plane, see \cite{BJP}. In addition, the work of Patterson-Perry \cite{PatPer} tells that the non-trivial zeros are the resonance set ${\mathcal R}_X$. The largest zero (in term of its real part) of $Z_\Lambda(s)$ is at $s=\delta$
and is known to be simple. There are no other zeros on the line $\{ \Re(s)=\delta \}$. To investigate zeros of Selberg zeta function
one can try to use the fact that it is a Fredholm determinant of a holomorphic family of trace class operators called "Ruelle transfer operators" $L_s$ so that 
$$Z_\Lambda(s)=\det(I-L_s).$$
Let us be more specific. Let $\D_1,\ldots,\D_p;\D_{p+1},\ldots \D_{2p}$, $p\geq 2$ be $2p$ closed, disjoints, euclidean discs
in $\C$ which we take to be orthogonal to $\partial \H$ (we use the Poincar\' e disc model). Let 
$h_1,\ldots,h_p \in \mathrm {PSU}(1,1)$ be such that for all $i=1,\ldots,p$
$$h_i(\D_i)=\overline{\widehat{\C}\setminus \D_{p+i}}.$$
The (free) group $\Lambda$ generated by $h_1,\ldots, h_p;h_1^{-1},\ldots,h_p^{-1}$ is convex co-compact 
(the converse is also true, every convex co-compact subgroup of $\mathrm {PSU}(1,1)$ can be realized as a Schottky group, see for example \cite{Borthwick} chapter 15). For convenience we set $h_{p+i}=h_{i}^{-1}$ for all $i=1,\ldots,p$. For all $i$, set $I_i:=\partial \H \cap \D_i$. The Bowen-Series map is an (eventually) expanding map defined on $\cup_{i=1}^{2p} I_i$ by $B(x)=h_i(x)$ if $x\in I_i$. The maximal invariant subset of $\partial \H$ for
the Bowen map $B$ is the limit set of $\Lambda$.
Set $$\U=\bigcup_{i=1}^{2p} \D_i,$$
and define the function space $\A(\U)$ to be 
$$\A(\U):=\{ f\in L^2(\U)\ :\ f\ \mathrm{is\ holomorphic\  on}\ \U \}.$$
The usual Ruelle transfer operator $L_s:\A(\U)\rightarrow \A(\U)$ is defined (see \cite{GLZ}) for $z\in \D_i$ by
$$ L_s(f)(z):=\sum_{j\neq i} (h_j')^s f(h_j z).$$
Note that all each $h_j:\D_i \rightarrow \D_{p+j}$ (indexes have to be understood mod $2p$)  is a contraction.
It has been known since Ruelle \cite{Ruelle0} that these operators are compact, trace class. The computation of the
trace yields
$$\mathrm{Tr}(L_s^n)=\sum_{\vert \alpha \vert =n} \frac{e^{-s l_\alpha}}{1-e^{-l_\alpha}},$$
where $2\cosh(l_\alpha/2)=\mathrm{Tr}(h_\alpha)$ and 
$h_\alpha=h_{\alpha_1}\ldots h_{\alpha_n} \in \Lambda$ has reduced word length $\vert \alpha \vert =n$
and $\alpha_n\neq \alpha_1+p \mod 2p$. Now if we write $s=\sigma+it$
$$\vert \mathrm{Tr}(L_s^n) \vert^2=\sum_{\vert \beta \vert,\ \vert \alpha \vert =n} 
\frac{e^{-\sigma (l_\alpha+l_\beta)-it(l_\alpha-l_\beta)}}{(1-e^{-l_\alpha})(1-e^{-l_\beta})},$$
and believe that the off-diagonal terms as $t\rightarrow +\infty$ will cancel up to the size of the diagonal sum,
we get 
$$\vert \mathrm{Tr}(L_s^n) \vert^2 =
O\left( \sum_{\vert \alpha \vert =n} \frac{e^{-2\sigma l_\alpha}}{(1-e^{-l_\alpha})^2}\right).$$
It known that this above sum is related to a quantity called "topological pressure" i.e. as $n\rightarrow \infty$,
$$\lim_{n\rightarrow +\infty} 
\left( \sum_{\vert \alpha \vert =n} \frac{e^{-2\sigma l_\alpha}}{(1-e^{-l_\alpha})^2}\right)^{1/n}=e^{P(-2\sigma)}.$$
The topological pressure is a thermodynamical quantity related to the Bowen map on the boundary by the following
formula, ($x\in [0,1]$)
$$P(-x)=\sup_{\mu} \left (h_\mu(B)-x\int_\Lambda \log\vert B' \vert d\mu\right),$$
where the sup is understood over all $B$-invariant measures and $h_\mu$ is the measure theoretic entropy. The map $x\mapsto P(-x)$ is strictly decreasing, convex, and 
Bowen's formula (see \cite{Bowen1}) says that it has a unique zero at $x=\delta$. Therefore it seems that
the series expansion
$$\det(I-L_s)=\exp \left(-\sum_{n=1}^{+\infty}\frac{1}{n} \mathrm{Tr}(L_s^n)\right)$$
is still converging for $\Re(s)>\delta/2$ and $\Im(s)$ large, which is in favor of $G(\Lambda)\leq \delta/2$.
On the other hand, it follows rigourously from \cite{NaudHeat} that for all $\Re(s)=\sigma$ fixed, 
for all $\epsilon>0$, one can find a sequence $(t_k)$ with $\vert t_k \vert \rightarrow +\infty$ as $k\rightarrow +\infty$
such that the spectral radius $\rho(L_{\sigma+it_k})$ satisfies the lower bound
$$\rho(L_{\sigma+it_k})\geq e^{3/2P(-2\sigma)},$$
which shows that for $\Re(s)<\delta/2$, the operators $L_s$ are no longer uniform contractions as 
$\vert \Im(s) \vert $ goes to infinity, which is in favor of $G(\Lambda)\geq \delta/2$.
\end{section}
\begin{section}{Proof of the main result}
We assume throughout this section that $\Lambda$ is a convex co-compact Fuchsian group which is non elementary. This causes no
trouble since the elementary case ($\delta=0$) is well understood in term of resonances \cite{ChristZworski}.
Although we could obtain Theorem \ref{mainresult} from the Guillop\'e-Zworski $0$-trace formula \cite{GuiZwor2}, it will prove more efficient for us to work directly with Selberg's zeta function instead. 
\subsection{An approximate trace formula}
Working with Selberg's zeta function plus an estimate due to
Guillop\'e and Zworski allows us to obtain an approximate trace formula. Let $\varphi \in C_0^\infty(\R)$ be a smooth compactly supported test function. We define $\psi(s)$ by the formula
$$\psi(s):=\int_{-\infty}^{+\infty}e^{su}\varphi(u)du=\widehat{\varphi}(is),$$
where $\widehat{\varphi}$ is the usual Fourier transform of $\varphi$. We have the following claim.
\begin{propo} 
\label{approx}
Assume that for a given $\rho<\delta$, the set 
$${\mathcal F}_\rho:={\mathcal R}_X\cap\{ \Re(s)> \rho \}$$ is finite,
then for all $\varepsilon>0$ small enough, one can find ${\varepsilon}\leq \widetilde{\varepsilon} \leq 2\varepsilon$ such that 
$$\sum_{k \in \N_0} \sum_{\gamma \in {\mathcal P}} \frac{l(\gamma)}{1-e^{-kl(\gamma)}}\varphi(kl(\gamma))=
\sum_{\lambda \in {\mathcal F}_\rho} \psi(\lambda)$$
$$+O\left ( \int_{-\infty}^{+\infty} (1+\vert x\vert)^\delta \vert \psi(\rho+\widetilde{\varepsilon}+ix)\vert dx \right).$$
The implied constant depending only on $\varepsilon$, $\rho$ and $\Lambda$.
\end{propo}
\noindent {\it Proof}. First we recall that for $\Re(s)>\delta$, we have the convergent series formula for the logarithmic derivative
$$\frac{Z_\Lambda'(s)}{Z_\Lambda(s)}=\sum_{ (n,\gamma)\in \N\times {\mathcal P}} \frac{l(\gamma)e^{-(s+n)l(\gamma)}}{1-e^{-(s+n)l(\gamma)}}=\sum_{k\in \N_0}\sum_{\gamma \in {\mathcal P}}\frac{l(\gamma)e^{-skl(\gamma)}}{1-e^{-kl(\gamma)}}.$$
In addition, for all $A>\delta$, one can find a constant $C_A\geq 0$ such that for all $\Re(s)\geq A$, we have
$$ \left \vert \frac{Z_\Lambda'(s)}{Z_\Lambda(s)} \right \vert\leq C_A.$$
Let us consider the contour integral
$$I_A:=\frac{1}{2i\pi}\int_{A-i\infty}^{A+i\infty} \frac{Z_\Lambda'(s)}{Z_\Lambda(s)}\psi(s)ds,$$
where $\psi(s)=\widehat{\varphi}(is)$. Since $\varphi \in C_0^\infty(\R)$, $\psi(s)$ is rapidly decreasing (Schwartz class)
on every vertical line which guarantees convergence of the above integral. Using Lebesgue dominated convergence theorem, we can
interchange $\sum$ and $\int$ to write 
$$I_A=\sum_{k\in \N_0}\sum_{\gamma \in {\mathcal P}} \frac{l(\gamma) e^{-Akl(\gamma)}}{1-e^{-kl(\gamma)}}
\times \frac{1}{2\pi}\int_{-\infty}^{+\infty} \widehat{e^{A x}\varphi(x)}(\xi) e^{i\xi kl(\gamma)}d\xi.$$
Fourier inversion formula yields
$$I_A=\sum_{k\in \N_0}\sum_{\gamma \in {\mathcal P}} \frac{l(\gamma)}{1-e^{-kl(\gamma)}} \varphi(kl(\gamma)).$$
Assuming that $Z_\Lambda(s)$ has only finitely many zeros inside $\{ \Re(s) >\rho \}$, the next step is to deform the contour integral
$I_A$. To perform this contour shift, we obviously need an a priori estimate on the growth of $Z_\Lambda'(s)/Z_\Lambda(s)$.
We fix $\varepsilon>0$. 
We show the following.
\begin{lem}
There exist a constant $C>0$ and $R_0\geq 0$, such that for all $\Re(s) \geq \rho+\varepsilon$
 and $\vert \Im(s) \vert \geq R_0$, we have 
 \begin{equation}
 \label{apriori1}
\left \vert \frac{Z_\Lambda'(s)}{Z_\Lambda(s)} \right \vert\leq C \vert \Im(s) \vert ^\delta.
\end{equation}
\end{lem}
\noindent {\it Proof}. From the scale-analysis of Guillop\'e-Lin-Zworski \cite{GLZ}, it is known that in every half-plane 
$\{ \Re(s)\geq \sigma \}$, one can find $C_\sigma\geq 0$ such that 
\begin{equation}
\label{fractal}
\vert Z_\Lambda(s)\vert \leq C_\sigma e^{C_\sigma \vert \Im (s) \vert ^\delta}.$$
For all $\vert t \vert$ large, the holomorphic function
$$z\mapsto \frac{Z_\Lambda(z+A+it)}{Z_\Lambda(A+it)}
\end{equation}
does not vanish on a neighborhood of the closed disc $\overline{D}(0,A-\rho-\varepsilon/2)$ so that we can define
a complex holomorphic logarithm
$$f_{t}(z):=\log\left ( \frac{Z_\Lambda(z+A+it)}{Z_\Lambda(A+it)} \right)$$
with $f_t(0)=0$ and 
$$\Re(f_t(z))=\log\left \vert \frac{Z_\Lambda(z+A+it)}{Z_\Lambda(A+it)} \right\vert.$$
We can now apply the Borel-Caratheodory estimate that can be found in Titchmarsh \cite{Tit}.
\begin{lem}
 Assume that $f$ is a holomorphic function on a neighborhood of the closed disc $\overline {D}(0,R)$, and $f(0)=0$. Then
 for all $r<R$, we have
 $$ \max_{\vert z\vert\leq r} \vert f'(z) \vert \leq \frac{8R}{(R-r)^2}\max_{\vert \zeta \vert \leq R} \vert \Re(f(z))\vert.$$
 \end{lem}
 \noindent
 Applying the above estimate to $f_t$, 
we eventually obtain for all $s$ in the closed disc $\overline{D}(A+it, A-\rho-\varepsilon)$, with $\vert t \vert$ large, 
$$\left \vert \frac{Z_\Lambda'(s)}{Z_\Lambda(s)} \right \vert
\leq \frac{16A}{\varepsilon}\left( C_\rho \vert \Im (s)\vert^\delta+\log C_\rho-\log \vert Z_\Lambda(A+it) \vert\right).$$
Since $\vert Z_\Lambda(s)\vert$ is uniformly bounded from below on $\{ \Re(s) \geq A \}$, we have proved 
the lemma. $\square$ 

\bigskip \noindent We can now conclude the proof of the approximate trace formula. 
We can choose $\widetilde{\varepsilon}\in [\varepsilon,2\varepsilon]$ such that 
$Z_\Lambda(s)$ does not vanish on the line $\{ \Re(s)=\rho+ \widetilde{\varepsilon}\}$.
Making sure that $\varepsilon$ is small enough such that 
$$\inf_{\lambda \in {\mathcal F}_\rho} \Re(\lambda)> \rho+2\varepsilon,$$ by the residue theorem we write for $R$ large,
$$ I_{A,R}:=\frac{1}{2i\pi}\int_{A-iR}^{A+iR} \frac{Z_\Lambda'(s)}{Z_\Lambda(s)}\psi(s)ds$$
$$=\sum_{\lambda \in {\mathcal F}_\rho} \psi(\lambda) + 
\underbrace{\frac{1}{2i\pi}\int_{\rho+\widetilde{\varepsilon}-iR}^{\rho+\widetilde{\varepsilon}+iR} \frac{Z_\Lambda'(s)}{Z_\Lambda(s)}\psi(s)ds}_{\widetilde{I}_{\rho,R}}
-\underbrace {\frac{1}{2i\pi}\int_{\Im(s)=\pm R \atop
\rho+\widetilde{\varepsilon} \leq \Re(s) \leq A
 }\frac{Z_\Lambda'(s)}{Z_\Lambda(s)}\psi(s)ds}_{E_R}.$$
Using the bound (\ref{apriori1}) plus the rapid decay of $\psi(s)$ on vertical lines, we get that $E_R\rightarrow 0$ as 
$R\rightarrow + \infty$. We therefore have the desired formula by taking the limit
as $R \rightarrow +\infty$ and use again the bound (\ref{apriori1}) in $\widetilde{I}_{\rho,R}$. $\square$ 

\subsection{Facts on the arithmetic length spectrum}
Given a discrete Fuchsian group $\Lambda$, one defines the lengh spectrum of $\Lambda\backslash \H$
as the set
$${\mathcal L}_\Lambda=\{ kl(\gamma)\ :\ (k,\gamma)\in \N_0\times {\mathcal P} \}.$$
If we view $\Lambda$ as a subgroup of ${\rm PSL}_2(\R)$, then there is a one-to-one correspondence between closed geodesics
on $X$ and conjugacy classes of hyperbolic elements in $\Lambda$. If $\gamma$ is a closed geodesic (not necessarily primitive)
corresponding to $g_\gamma \in \Lambda$, a representative of a conjugacy class, then its length $l(\gamma)$ is related to
the trace of $g_\gamma$ through the formula
\begin{equation}
\label{trace1}
\vert {\rm Tr}(g_\gamma) \vert=2\cosh(l(\gamma)/2).
\end{equation}
We now recall a few basic facts on the properties of the length spectrum of arithmetic groups. An arithmetic Fuchsian group
is often defined as a discrete subgroup of ${\rm PSL}_2(\R)$ which is commensurable to a subgroup derived from a quaternion
algebra, we refer to the book of Katok \cite{Katok} for precise definitions. For our purpose, we will only need the following
characterization due to Takeuchi \cite{Takeuchi}. 
\begin{thm}[Takeuchi]
Let $\Lambda$ be a discrete, cofinite subgroup of $\textrm{PSL}_2(\R)$. Set
${\rm Tr}(\Lambda):=\{ {\rm Tr}(T)\ :\ T\in \Lambda \}$, ${\rm Tr}(\Lambda)^2:=\{ ({\rm Tr}(T))^2\ :\ T\in \Lambda \}$.
Then $\Lambda$ is arithmetic if and only if:
\begin{enumerate}
 \item The field $K=\Q( {\rm Tr}(\Lambda) )$ is an algebraic field of finite degree and ${\rm Tr}(\Lambda)$ is a subset of the ring of integers
 of $K$.
 \item For all embedding $\varphi:K\rightarrow \C$ with $\varphi \vert_{{\rm Tr}(\Lambda)^2} \neq Id$, the set $\varphi( {\rm Tr}(\Lambda))$ is bounded in $\C$.
\end{enumerate} 
\end{thm}
\noindent
This special features of the trace set ${\rm Tr}(\Lambda)$ imply the important {\it bounded cluster} property: there exists a constant
$C_\Lambda>0$ such that for all $n\in \N$,
$$ \# {\rm Tr}(\Lambda)\cap [n,n+1]\leq C_\Lambda.$$ A proof of that fact (based on the number theoretic properties above) can be found in the paper of Luo and Sarnak \cite{SarnakLuo}. The main consequence for our purpose is the following lower bound.
Let $\Lambda$ be a convex co-compact subgroup of an arithmetic group $\Lambda_0$.
For all $\ell \in {\mathcal L}_\Lambda$, let $m(\ell)$ denote the multiplicity of $\ell$ as the length of a closed geodesic i.e.
$$m(\ell):=\#\{ (k,\gamma)\in \N_0\times {\mathcal P}\ :\ \ell=kl(\gamma) \}.$$
We have the next result.
\begin{propo}
\label{mult} 
Assume that $\delta(\Lambda) >\half$, then
there exists $A_\Lambda>0$ such that for all $T$ large, we have
$$\sum_{T-1\leq \ell \leq T+1 \atop \ell \in {\mathcal L}_\Lambda}
 m^2(\ell)\geq A_\Lambda \frac{e^{(2\delta(\Lambda)-1/2)T}}{T^2}.$$
\end{propo}
\noindent {\it Proof}. Because of the bounded cluster property, we have
$$\# \{ t \in {\rm Tr}(\Lambda)\ :\ t\leq x\}\leq C_{\Lambda_0} [x],$$
which using formula (\ref{trace1}) gives
$$\# \{ \ell \in {\mathcal L}_\Lambda\ :\ \ell \leq T\}\leq 2C_{\Lambda_0}\cosh(T/2)\leq C' e^{T/2},$$
for all $T$ large and a well chosen $C'$.
On the other hand, the prime geodesic theorem (which is due to Lalley \cite{Lalley} for convex co-compact groups) 
applied to $X=\Lambda\backslash \H$ states that 
$$\#\{ (k,\gamma) \in \N_0\times {\mathcal P}\ :\ kl(\gamma)\leq T \}=\frac{e^{\delta T}}{\delta T}\left( 1+o(1)\right),$$
which gives
$$\sum_{T-1\leq \ell \leq T+1 \atop \ell \in {\mathcal L}_\Lambda} m(\ell)\geq \widetilde{C} \frac{e^{\delta T}}{ T},$$
for some constant $\widetilde{C}>0$ depending on $\Lambda$. Applying Schwarz inequality gives
$$\widetilde{C} \frac{e^{\delta T}}{T}\leq 
\sqrt{C'}e^{T/4} \left (\sum_{T-1\leq \ell \leq T+1 \atop \ell \in {\mathcal L}_\Lambda} m^2(\ell) \right)^{1/2} ,$$
and the proof is finished. $\square$
\subsection{Main proof}
We can now proceed to the proof of Theorem \ref{mainresult}. The idea is to work by contradiction and use the approximate
trace formula of Proposition \ref{approx}. We will use the following family of test functions with parameters $\xi \in \R$ and
$T\geq 0$, $T$ large.
$$\varphi_{\xi,T}(x):=e^{-i\xi x}\varphi(x-T),$$
where $\varphi \in C_0^\infty(\R)$ with support inside $[-2,2]$, $\varphi \geq 0$ and $\varphi(x)=1$ for all $x\in [-1,1]$.
We set for $A\geq \Re(s)\geq 0$,
$$\psi_{\xi,T}(s):=\widehat{\varphi_{\xi,T}}(is)=e^{-i\xi T}e^{sT}\widehat{\varphi}(\xi+is).$$
Notice that by integrating by parts, for all integer $N\geq 0$, one can find $C_N>0$ such that
\begin{equation}
\label{rapid}
\vert \psi_{\xi,T}(s) \vert \leq C_N \frac{e^{\Re( s)T}}{(1+\vert \xi +\Im(s) \vert)^N}. 
\end{equation}
Assume from now on that there exist only finitely many resonances inside the strip $\{ \Re(s)>\rho\}$.  We set  
$${\mathcal S}_{\xi,T}:=\sum_{k \in \N_0} \sum_{\gamma \in {\mathcal P}} \frac{l(\gamma)}{1-e^{-kl(\gamma)}}
e^{-i\xi kl(\gamma)}\varphi(kl(\gamma)-T).$$
Applying the approximate trace formula to $\varphi_{\xi,T}$ we have
\begin{equation}
\label{eq1}
{\mathcal S}_{\xi,T=}\sum_{\lambda \in {\mathcal F}_\rho} \psi_{\xi,T}(\lambda)
+E(\xi,T) ,
\end{equation}
where the error term $E(\xi,T)$ satisfies
$$E(\xi,T)=O\left (\int_{-\infty}^{+\infty} (1+\vert x\vert)^\delta \vert \psi_{\xi,T}(\rho+\widetilde{\varepsilon}+ix)\vert dx\right).$$
Using the estimate (\ref{rapid}) with $N=2$ and the fact that $\delta<1$, we can bound this error term by 
$$E(\xi,T)=O\left (  e^{(\rho+\widetilde{\varepsilon})T} I(\xi) \right),$$
where $I(\xi)$ is the following convergent integral
$$I(\xi)=\int_{-\infty}^{+\infty} \frac{(1+\vert x \vert )^\delta}{(1+\vert \xi +x\vert)^2 }dx.$$
An obvious change of variable shows that for $\vert \xi \vert \geq 1$, $I(\xi)$ is bounded from above by
$$ I(\xi)=\int_{-\infty}^{+\infty} \frac{(1+\vert u-\xi \vert )^\delta}{(1+\vert u\vert)^2 }du\leq \vert \xi \vert ^\delta
\int_{-\infty}^{+\infty} \frac{(2+\vert u\vert )^\delta}{(1+\vert u\vert)^2 }du,$$
so that we have
\begin{equation}
\label{eq2}
E(\xi,T)=O\left (  e^{(\rho+\widetilde{\varepsilon})T} \vert \xi \vert^\delta \right).
\end{equation}
The key step is now to consider the weighted squared sum 
\footnote{There is nothing sacred about the gaussian weight here. One can use any "squared convolued" weight $W_\sigma(\xi)$
of the form $$W_\sigma(\xi)=f \star f(\sqrt{\sigma}\xi),$$ where $f\geq 0$ is in Schwartz class so that the Fourier transform
of $W_\sigma(\xi)$ is positive.}
$${\mathcal G}(\sigma,T):=\sqrt{\sigma}\int_{-\infty}^{+\infty} e^{-\sigma\xi^2} \vert {\mathcal S}_{\xi,T} \vert^2 d\xi,$$
where $\sigma>0$ is a small parameter depending on $T$ (which is large) to be adjusted later on. 
We recall that ${\mathcal L}_\Lambda$ stands for the length spectrum. If $\ell \in {\mathcal L}_\Lambda$, we denote
$\widetilde{\ell}$ the {\it prime length} of $\ell$, i.e. if $\ell=kl(\gamma)$ with $\gamma \in {\mathcal P}$, then 
$\widetilde{\ell}=l(\gamma)$. We also set for $\ell, \ell' \in {\mathcal L}_\Lambda$, 
$$a_{\ell,\ell'}:=\frac{\widetilde{\ell} \widetilde{ \ell'}m(\ell)m(\ell')}{(1-e^{-\ell})(1-e^{-\ell'})}.$$
Interchanging summations we have
$${\mathcal G}(\sigma,T)=\sum_{\ell,\ell' \in {\mathcal L}_\Lambda} a_{\ell,\ell'}\varphi(\ell-T)\varphi(\ell'-T)\int_{-\infty}^{+\infty}
e^{-i\xi(\ell-\ell') }e^{-\sigma\xi^2}\sqrt{\sigma}d\xi,$$
using the classical formula for the Fourier transform of Gaussian functions, we obtain
$${\mathcal G}(\sigma,T)=\sqrt{\pi} \sum_{\ell,\ell' \in {\mathcal L}_\Lambda} a_{\ell,\ell'}\varphi(\ell-T)\varphi(\ell'-T) 
e^{-\frac{(\ell-\ell')^2}{4\sigma}}.$$ 
Since all summands in the above formula are positive, we can simply drop all off-diagonal terms to get the lower bound
$${\mathcal G}(\sigma,T)\geq \sqrt{\pi} \sum_{\ell \in {\mathcal L}_\Lambda} a_{\ell,\ell}\varphi^2(\ell-T),$$
which eventually yields
$${\mathcal G}(\sigma,T)\geq C \sum_{T-1\leq \ell \leq T+1 \atop \ell \in {\mathcal L}_\Lambda}
 m^2(\ell),$$
 for all $T$ large and a suitable constant $C>0$.
 Going back to (\ref{eq1}), we have 
 $$\frac{{\mathcal G}(\sigma,T)}{2\sqrt{\sigma}}\leq \underbrace{\int_{-\infty}^{+\infty}e^{-\sigma\xi^2} \left \vert \sum_{\lambda \in {\mathcal F}_\rho} 
 \psi_{\xi,T}(\lambda)\right \vert^2 d\xi}_{{\mathcal I}_1(\sigma,T)} +\underbrace{\int_{-\infty}^{+\infty}e^{-\sigma\xi^2} \vert E(\xi,T)\vert^2 d\xi}
 _{{\mathcal I}_2(\sigma,T)}.$$
 Using the estimate (\ref{eq2}), one obtains that 
 $${\mathcal I}_2(\sigma,T)=O\left ( e^{2(\rho+\widetilde{\varepsilon})T}\int_{-\infty}^{+\infty} e^{-\sigma\xi^2}\vert \xi \vert^{2\delta} d\xi\right)$$
 $$=O\left(  e^{2(\rho+\widetilde{\varepsilon})T}\sigma^{-\delta-\half}  \right).$$
 On the other hand, since ${\mathcal F}_\rho$ is finite and using estimate (\ref{rapid}) with $N=1$, one can find 
 $\widetilde{C}>0$ such that for all $\lambda \in {\mathcal F}_\rho$, 
 $$\vert \psi_{\xi,T}(\lambda) \vert \leq \widetilde{C} \frac{e^{\delta T}}{(1+\vert \xi \vert )},$$
 which shows that
 $${\mathcal I}_1(\sigma,T)=O\left ( e^{2\delta T}\int_{-\infty}^{+\infty}\frac{e^{-\sigma\xi^2}}{(1+\vert \xi\vert)^2}d\xi \right)=
 O\left ( e^{2\delta T} \right),$$
 uniformly in $\sigma$. We can now conclude the proof, depending on $\delta$.
 
 \bigskip \noindent
 If $0<\delta\leq \half$, we cannot use Proposition \ref{mult} and we have to rely on the weaker lower bound
 $${\mathcal G}(h,T)\geq C \sum_{T-1\leq \ell \leq T+1 \atop \ell \in {\mathcal L}_\Lambda}
 m(\ell),$$
 which by the prime orbit theorem yields
 $${\mathcal G}(\sigma,T)\geq B \frac{e^{\delta T}}{T},$$
 for some well chosen $B>0$. We have therefore obtained for all $T$ large and $\sigma$ small,
 $$B \frac{e^{\delta T}}{T}=O\left ( \sqrt{\sigma}e^{2\delta T} \right)+ O\left(  e^{2(\rho+\widetilde{\varepsilon})T}\sigma^{-\delta}  \right).$$
 Setting $\sigma=e^{-\alpha T}$, we have a contradiction as $T\rightarrow +\infty$ provided
 $$\alpha>2\delta\ {\rm and}\ \rho< \frac{\delta(1-\alpha)}{2}-\widetilde{\varepsilon},$$
 which shows ($\widetilde{\varepsilon}$ being as small as one wants) that there are infinitely many resonances in the strip $\{\Re(s)\geq \frac{\delta(1-2\delta)}{2}-\epsilon \}$ for all $\epsilon>0$.
 
\bigskip \noindent
 If $\delta>\half$ , we use the lower bound of Proposition \ref{mult} to write (for some $B>0$)
 $$B \frac{e^{(2\delta-\half) T}}{T^2}=O\left ( \sqrt{\sigma}e^{2\delta T} \right)+ O\left(  e^{2(\rho+\widetilde{\varepsilon})T}
 \sigma^{-\delta}  \right),$$
 which if $\sigma=e^{-\alpha T}$ produces a contradiction whenever $\alpha>1$ and
 $$ \rho<\frac{\delta(2-\alpha)}{2}-\frac{1}{4}-\widetilde{\varepsilon}.$$
 As above, we have found infinitely many resonances (or zeros of $Z_\Lambda(s)$) in the strip
 $\{\Re(s)\geq \frac{\delta}{2}-\frac{1}{4}-\epsilon \}$ for all $\epsilon>0$.  $\square$
  
\bigskip
Finally we would like to point out an alternative proof ($\delta >1/2$) if $\Lambda$ is an  infinite index subgroup of an arithmetic
group $\Lambda_0$ derived from quaternion algebra. In that specific case, a stronger version of the bounded cluster property is valid: there exist $C_{\Lambda_0}>0$ such that given $g_1, g_2 \in \Lambda_0$ with $\mathrm{Tr}(g_1)\neq \mathrm{Tr}(g_2)$, 
$$\vert \mathrm{Tr}(g_1)-\mathrm{Tr}(g_2)\vert \geq C_{\Lambda_0}.$$
This property allows do deal with oscillatory terms. Consider the same family of test functions as above
and define instead
 $${\mathcal H}(R,T):=\frac{1}{R}\int_{R}^{3R}\left(1-\frac{\vert \xi-2R\vert}{R}\right)
 \vert{\mathcal S}_{\xi,T} \vert^2d\xi.$$
 The nice thing about these quantities is that 
 $$\frac{1}{R}\int_{R}^{3R}\left(1-\frac{\vert \xi-2R\vert}{R}\right)\left \vert \sum_{\lambda \in {\mathcal F}_\rho} 
 \psi_{\xi,T}(\lambda)\right \vert^2 d\xi=O\left( \frac{e^{2\delta T}}{R^{-\infty}} \right),$$
and
 $$\frac{1}{R}\int_{R}^{3R}\left(1-\frac{\vert \xi-2R\vert}{R}\right)\left \vert E(\xi,T)\right \vert^2 d\xi=
 O\left( e^{2(\rho+\widetilde{\varepsilon} )T}R^{2\delta} \right).$$
 However the price to pay for that is the formula
 $${\mathcal H}(R,T)=R^{-1}\sum_{\ell,\ell' \in {\mathcal L}_\Lambda} a_{\ell,\ell'}\varphi(\ell-T)\varphi(\ell'-T)
 G(R,\ell-\ell'),$$
 with $$G(R,x):=\int_{R}^{3R}\left( 1-\frac{\vert \xi-2R\vert}{R}\right)e^{-i\xi x }d\xi,$$
 which is no longer a sum of positive terms. Integrating by parts we found that the off-diagonal sum is of size at 
 most
$$R^{-2}\sum_{\ell,\ell' \in {\mathcal L}_\Lambda} 
a_{\ell,\ell'}\frac{\varphi(\ell-T)\varphi(\ell'-T)}{\vert \ell-\ell' \vert^2}.$$ 
In addition to that we know that 
$$\vert \ell-\ell' \vert \geq e^{-\max(\ell,\ell')/2}\vert  \mathrm{Tr}(g_\ell) -\mathrm{Tr}(g_{\ell'})\vert,$$ 
where $g_\ell \in \Lambda$ satisfies $2\cosh(\ell/2)=\mathrm{Tr}(g_\ell)$.
The diagonal contribution being proportional to 
$$\sum_{\ell \in {\mathcal L}_\Lambda} 
a_{\ell,\ell}\varphi^2(\ell-T),$$
we find that the off-diagonal contribution does not interfere 
with the diagonal one as long as $R\geq C e^{T/2}$,
for some well chosen constant $C>0$. In the end we are left with
$$\frac{e^{(2\delta-\half)T}}{T^2}=
O\left(  e^{(2(\rho+\widetilde{\varepsilon})+\delta )T}   \right), $$
which yields a contradiction for $\rho< \delta/2-1/4$.
 
\end{section}

\section{The error term in the hyperbolic lattice counting problem}
In this last section, we prove the omega lower bound for the error 
term in Lax-Phillips lattice counting asymptotic.
Let us assume that we have for fixed $z,z'$ in $\H$ with $z\neq z' $ 
in $\Lambda \backslash \H$ and $0<\rho <\delta$,
$$N(T;z,z')=\sum_j Q_j(T;z,z')e^{\delta_j T}+O\left(  e^{\rho T}\right),$$
where $j$ runs over a finite set, $\delta_j \in \C$ 
and $Q_j(T,z,z')=\sum_{k=0}^{N(j)} a_{k,j}(z,z')T^k$.
The usual Poincar\'e series defined for $\Re(s)>\delta$ by
$$P_\Lambda(s;z,z'):=\sum_{\gamma \in \Lambda} e^{-s {\mathrm d}(z,\gamma z')},$$
are related to the counting function $N(T;z,z')$ via summation by parts 
through the formula
$$P_\Lambda(s;z,z')=s\int_0^{+\infty}e^{-st}N(t;z,z')dt.$$
Inserting the asymptotic in this formula we get for $\Re(s)>\delta$,
$$P_\Lambda(s;z,z')= \sum_{j}\sum_{k=0}^{N(j)} 
a_{k,j}(z,z') \frac{sk!}{(s-\delta_j)^{k+1}}+F(s;z,z'),$$
where $s\mapsto F(s;z,z')$ has a holomorphic extension to the half plane 
$$\{ \Re(s)>\rho \}.$$
Poincar\'e series are in turn related to the resolvent kernel $R_\Lambda(s;z,z')$ on $\Lambda \backslash \H$ by the Gauss hypergeometric function
averaged over the group $\Lambda$. More precisely, for $\Re(s)>\delta$ and $z\neq z' $ in $\Lambda \backslash \H$, we have
$$R_\Lambda(s;z,z')=\sum_{\gamma \in \Lambda} G_s(\sigma(z,\gamma z')),$$
with $\sigma(z,z')=\cosh^2({\mathrm d}(z,z')/2)$ and for $\sigma>1$,
$$G_s(\sigma)=\frac{1}{4\pi}\sum_{k=0}^\infty \frac{(\Gamma(s+k))^2}{k!\Gamma(2s+k)}\sigma^{-s-k}.$$
>From that fact one infers that there exists a holomorphic function on the half plane
$$\{ \Re(s)>\delta-1\}$$
 denoted by $H(s;z,z')$ such that
$$R_\Lambda(s;z,z')=\frac{4^{s-1}}{\pi} \frac{(\Gamma(s))^2}{\Gamma(2s)}P_\Lambda(s;z,z')+H(s;z,z'),$$
see \cite{Borthwick}, proposition 14.17 for computational details.
We have therefore shown that for $(z,z')$ away from the diagonal in $\Lambda \backslash \H \times \Lambda \backslash \H$, the map
$$s\mapsto R_\Lambda (s;z,z')$$ extends meromorphically to $\{ \Re(s) >\rho \}$ with only finitely many poles.

On the other hand, we know from the parametrix of Guillop\'e-Zworski that (again away from the diagonal) 
$s\mapsto R_\Lambda(s;z,z')$ has a meromorphic extension to $\C$. The structure of this kernel at a resonance is as follows, see
\cite{Borthwick} proposition 8.1.
Assume that $\zeta\in \C$ is a pole, then there exist $p\geq 1$ and a finite family $A_j(\zeta;z,z')$ of smooth, finite rank kernels such that
$$R_\Lambda(s;z,z')=\sum_{j=1}^p\frac{A_j(\zeta;z,z')}{(s(1-s)-\zeta(1-\zeta))^j}+H(s;z,z'),$$
where $s\mapsto H(s;z,z')$ is holomorphic on a neighborhood of $\zeta$. The residue $A_1(\zeta;z,z')$ is always non-trivial and
can be described more
precisely as 
$$A_1(\zeta;z,z')=\sum_{k=1}^{m_\zeta}\phi_{k,\zeta}(z)\phi_{k,\zeta}(z'),$$
where each $\phi_{k,\zeta} \in C^\infty(X)$, $X=\Lambda\backslash \H$, is a non-trivial solution of the elliptic equation
$$(\Delta_X-s(1-s))^p \phi=0.$$ The coefficients of the elliptic operator $(\Delta_X-s(1-s))^p$ being obviously real-analytic,
by a well known theorem, each $\phi_{k,\zeta}$ is a (possibly $\C$-valued) non-trivial real-analytic function on $X$. Therefore
the function 
$$(z,z')\mapsto A_1(\zeta;z,z')$$
is real-analytic (and non trivial)
therefore its zero set ${\mathcal Z}(\zeta)$ has zero Lebesgue measure. Now fix an $\epsilon_0>0$ and consider  the set
$$ {\mathcal N}:=\bigcup_{\Re(\zeta)\geq \frac{\delta}{2}-\frac{1}{4}-\epsilon_0} 
{\mathcal Z}(\zeta) ,$$
where the union is understood over the (infinite) countable set of all resonances in the half-plane
$\{\Re(\zeta)\geq \frac{\delta}{2}-\frac{1}{4}-\epsilon_0\}$. Clearly the Lebesgue measure of ${\mathcal N}\subset X\times X$ is zero. 
Consider now 
$${\mathcal G}:=X\times X\setminus ({\mathcal N}\cup \Delta),$$
where $\Delta$ stands for the diagonal in $X\times X$. The set $\mathcal G$ is of full Lebesgue measure and for each
$(z,z')\in {\mathcal G}$ the kernel resolvent $s\mapsto R_\Lambda(s;z,z')$ has infinitely many poles in 
$\{\Re(s)\geq \frac{\delta}{2}-\frac{1}{4}-\epsilon\}$ for all $0<\epsilon\leq \epsilon_0$. We have a contradiction whenever 
$\rho\leq\frac{\delta}{2}-\frac{1}{4}-\epsilon$.

\end{document}